\documentclass[12pt,a4paper,twoside,english,final]{article}

\usepackage{babel}

\usepackage[utf8]{inputenc} 

\usepackage[T1]{fontenc}

\usepackage{microtype} 

\usepackage{lmodern}
\usepackage{ae}

\usepackage{graphicx} 

\usepackage{color}

\usepackage{float}

\usepackage{hypbmsec}

\usepackage[pdftex,hyperindex,colorlinks,bookmarks=false]{hyperref}
 
\usepackage{fancyhdr}

\usepackage{makeidx}

\usepackage[numindex,nottoc,numbib,section]{tocbibind} 

\usepackage[titles]{tocloft}    

\usepackage{showkeys}

\usepackage{changebar}

\usepackage{listings}
\definecolor{lgray}{gray}{0.80}
\lstnewenvironment{mupad}{\lstset{backgroundcolor=\color{lgray}, frame=single, basicstyle=\ttfamily, breaklines=true,mathescape =true}}{}

\usepackage[leqno]{amsmath}	

\usepackage{amssymb}		

\usepackage{amsxtra}		

\usepackage[thmmarks,hyperref]{ntheorem}

\usepackage{bm}                 

\usepackage{enumerate}
 
\usepackage{array}
\usepackage{eqnarray}

\usepackage{longtable}

\usepackage{diagrams}

\fancypagestyle{plain}{%
 \fancyhf{}
 }

\fancypagestyle{Index}{%
 \fancyhf{}%
 \fancyhead[LE,RO]{\bfseries\thepage}%
 \fancyhead[LO]{\bfseries Index}%
 \fancyhead[RE]{\bfseries Index}}

\fancypagestyle{Vorwort}{%
 \fancyhf{}%
 \fancyhead[LE,RO]{\bfseries\thepage}%
 \fancyhead[LO]{\bfseries Vorwort}%
 \fancyhead[RE]{\bfseries Vorwort}}

\fancypagestyle{Abbildungsverzeichnis}{%
 \fancyhf{}%
 \fancyhead[LE,RO]{\bfseries\thepage}%
 \fancyhead[LO]{\bfseries Abbildungsverzeichnis}%
 \fancyhead[RE]{\bfseries Abbildungsverzeichnis}}

\fancypagestyle{Inhaltsverzeichnis}{%
 \fancyhf{}%
 \fancyhead[LE,RO]{\bfseries\thepage}%
 \fancyhead[LO]{\bfseries Inhaltsverzeichnis}%
 \fancyhead[RE]{\bfseries Inhaltsverzeichnis}}

\fancypagestyle{Bibliographie}{%
 \fancyhf{}%
 \fancyhead[LE,RO]{\bfseries\thepage}%
 \fancyhead[LO]{\bfseries Literatur}%
 \fancyhead[RE]{\bfseries Literatur}}

\fancypagestyle{normal}{%
 \fancyhf{}%
 \fancyhead[LE,RO]{\bfseries\thepage}%
 \fancyhead[CE]{\bfseries Stefan Wiedmann}%
 \fancyhead[CO]{\bfseries Level-Structers of Drinfeld-Modules -- Closing a small gap}}

\theoremstyle{break}
\newtheorem{Theorem}{Theorem}[section]
\newtheorem{Satz}[Theorem]{Satz}
\newtheorem{Proposition}[Theorem]{Proposition}
\newtheorem{Lemma}[Theorem]{Lemma}

\newtheorem{Definition}[Theorem]{Definition}

\newtheorem{Corollary}[Theorem]{Corollary}

\theorembodyfont{\normalfont}

\newtheorem{Remark}[Theorem]{Remark}

\newtheorem{Properties}[Theorem]{Properties}
\newtheorem{Results}[Theorem]{Results}

\theoremstyle{nonumberplain}
\theorembodyfont{\normalfont}
\theoremsymbol{\ensuremath{_\Box}}
\newtheorem{Beweis}{Proof}

\newcommand{\OS}{\mathcal{O}}           
\newcommand{\Garbe}[1]{\mathcal{#1}}    

\newcommand{\field}[1]{\mathbb{#1}}     

\newcommand{\C}{\field{C}}              
\newcommand{\Fq}{\field{F}_{\!q}}       
\newcommand{\Z}{\field{Z}}              

\newcommand{\fdg}{\;|\;}                 

\newcommand{\prim}[1]{\mathfrak{#1}}   

\newcommand{\Gas}[1]{\field{G}_{a/#1}} 

\DeclareMathOperator{\AddPol}{AddPol}

\DeclareMathOperator{\Char}{char}

\DeclareMathOperator{\End}{End}

\DeclareMathOperator{\Hom}{Hom}

\DeclareMathOperator{\Kern}{Ker}

\DeclareMathOperator{\Quot}{Quot}

\DeclareMathOperator{\rk}{rk}

\DeclareMathOperator{\spec}{spec}

\usepackage{parskip}

\title{Level-Structures of Drinfeld-Modules -- Closing a small gap\\[2cm]}
\author{Stefan Wiedmann}
\date{Göttingen 2008}

\begin{document}


\maketitle
\tableofcontents          

\pagestyle{normal}      
\section*{Introduction}
Level structures play an important role in the definition of moduli spaces, because they offer a possibility to rigidify moduli problems. For example level $N$ structures of elliptic curves over $\C$ are defined by an isomorphism of $\Z/N \times \Z/N$ and $N$-torsion points. If one replaces the base $\C$ by an arbitary scheme then one is forced to regard the $N$-torsion points as a group scheme which possibly has at some points of the base connected components. So in this case the concept of an isomorphism of a constant group scheme and a the group scheme of $N$-division points does not work any more. After Drinfeld one weakens this isomorphism condition to a morphism which matches the corresponding Cartier divisors properly. This idea leads to the notion of \emph{generators} of a level structure \cite[(3.1)]{KatzMazur}, respectively in a more general setup to the notion of \emph{$A$-structures} and \emph{$A$-generators} \cite[Ch.~1]{KatzMazur}.

A similar situation occurs in the theory of Drinfeld modules, which can be seen as an analogue of elliptic curves in characteristic $p$. As it is discussed in loc.~cit.~one runs into a small difficulty in the definition of level structures as a condition for all prime divisors of $N$ versus a condition for $N$ itself which is discussed in loc.~cit.~in general \cite[Prop.~1.11.3, Rk.~1.11.4]{KatzMazur} and in the case of elliptic curves \cite[Th.~5.5.7]{KatzMazur}. The aim of this article is to prove the analogous result in the case of Drinfeld modules.

An  $I$-level structure of a Drinfeld module $(E,e)$ of rank $d$ over a base scheme $S$ is defined as an $A$ module homomorphism
\[
 \iota : \left(I^{-1}/A\right)^d \rTo E(S)
\]
such that for all prime ideals $\mathfrak{p} \supseteq I$ we have an equaltity
\[
 E[\mathfrak{p}] = \sum_{x \in \left(\mathfrak{p}^{-1}/A\right)^d} \iota(x)
\]
 of relative Cartier divisors \cite{drinfeld76}. In \cite{lehmkuhl} 4, prop. 3.3, it is proved that if $\iota$ 
 is an $I$-level structure then
 \[
 E[I] = \sum_{x \in \left(I^{-1}/A\right)^d} \iota(x)
\]
 is an equality of relative Cartier divisors too. 

 If one defines an $I$-level structure by the equality 
 \[
 E[I] = \sum_{x \in \left(I^{-1}/A\right)^d} \iota(x)
\]
of relative Cartier divisors does it follow that
\[
 E[\mathfrak{p}] = \sum_{x \in \left(\mathfrak{p}^{-1}/A\right)^d} \iota(x)
\]
is an equality of relativ Cartier divisors for every prime ideal $\mathfrak{p} \supseteq I$ too?. This question was solved in the authors PhD-Thesis \cite{wiedmanndiss} if the base scheme $S$ is reduced. A careful reading of the arguments used in \cite{lehmkuhl} gives the general result for an arbitary base scheme $S$.

\section{Definitions}
\subsection{Drinfeld modules}
Let $X$ be a geometrically connected smooth algebraic curve over the finite field
$\Fq$, let $\infty \in X$ be a closed point and let
$A := \Gamma(X \setminus \infty, \OS_X)$ be the ring of regular functions
outside $\infty$. In this case $A$ is
a 
Dedekind ring. 

Let $S/\Fq$ be a scheme, $\Garbe{L}$ a line bundle over $S$
and let $\Gas{\Garbe{L}}$ be the additive group scheme corresponding to the line
 bundle
$\Garbe{L}$. For all open subsets $U \subset S$ the 
group scheme is defined by
\begin{displaymath}
  \Gas{\Garbe{L}}(U) = \Garbe{L}(U).
\end{displaymath}
The (additive) groups $\Garbe{L}(U)$ are in a canonical way 
$\Fq$-vector spaces.

\begin{Definition}[\cite{drinfeld76}]\label{1:1}
  Let $\Char: S \rTo \spec A$ be a morphism over $\Fq$. A
  \emph{Drinfeld module} $E:= (\Gas{\Garbe{L}},e)$ consists of an 
    additive
  group scheme $\Gas{\Garbe{L}}$ and a ring homomorphism
   \[
   e: A \rTo \End_{\Fq}(\Gas{\Garbe{L}})
   \]
   such that:
  \begin{enumerate}[1)]
  \item The morphism $e(a)$ is finite for all $a \in A$ and for
    all points $ s \in S$ there exists an element $a \in A$ such that locally in $s$ the rank of the morphism $e(a)$ is 
    bigger than $1$.
  \item The diagram
    \begin{diagram}[silent]
      A&          &\rTo^e&          & \End_{\Fq}(\Gas{\Garbe{L}})\\
      &\rdTo_{\Char}&      &\ldTo_{\partial}&\\
      &&\OS_S(S)&&
    \end{diagram}
    commutes.
  \end{enumerate}
If $S=\spec R$ is affine und if $\Garbe{L}$ is trivial, then we will simply write $E = (R,e)$. 
\end{Definition}

\begin{Proposition}
Let $S$ be a connected scheme. Then there exists a natural number $d >0$, such that for all $0 \neq a \in A$ we have $\rk(e(a)) = q^{-d\deg (\infty)\infty(a)}$. The number $d$ is called the \emph{rank} of the Drinfeld module.
\end{Proposition}

If $S = \spec(R)$ is an affine scheme and if $\Garbe{L}$ is trivial, then we can show that 
\[
\End_{\Fq}(\Gas{\Garbe{L}}) \cong \AddPol_q(R)
\]
where $\AddPol_q(R)$ is the ring of $\Fq$-linear polynomials, i.e.~every polynomial $f(X) \in \AddPol_q(R)$ is of the form
\[
f(X) = \sum_{i=0}^n \lambda_i X^{q^i}.
\] 
In the affine situation a Drinfeld module is therefore given by a non trivial ring homomorphism $e$ and a commutative diagram
\begin{diagram}[silent]
      A&          &\rTo^e&          & \AddPol_q(R)\\
      &\rdTo_{\Char}&      &\ldTo_{\partial}&\\
      &&R&&
\end{diagram}
If we define $e_a(X) := e(a)$ and if $e_a(X) = \sum_{i=0}^n \lambda_i X^{q^i}$ then $\partial(e_a(X)) = \lambda_0$, the coefficient $\lambda_{-d\deg (\infty)\infty(a)}$ is a unit in $R$ and $\lambda_i$ is nilpotent for $i > -d\deg (\infty)\infty(a)$. If in this case $\lambda_i = 0$ for all $a \in A$ then the Drinfeld module is called \emph{standard}. One can show, that every Drinfeld module is isomorphic to a Drinfeld module in standard form.

By abuse of language the image of the map $\Char : S \rTo \spec A$ is called the \emph{characteristic} of $E$. If it consists only of the zero ideal then we say $E$ has \emph{general characteristic}.

\subsection{Division points and level structures}
Let $E$ be a Drinfeld module of rank $d$ over a base scheme $S$ and let $0 \neq I \subsetneq A$ be an ideal. 
\begin{Definition}
  The contravariant functor $E[I]$ on the category of schemes 
  over $S$ with image in the category of $A/I$ modules defined by
  \begin{displaymath}
    T/S \rMapsto \{ x \in E(T) \fdg I x = 0\} = \Hom_A(A/I, E(T))
  \end{displaymath}
  for all schemes $T/S$ is called the scheme of $I$-division points. 
\end{Definition}

\begin{Properties}
  \begin{enumerate}[1)]
  \item $E[I] \subseteq E$ is a closed, finite and flat (sub-)group scheme over $S$ of rank $|A/I|^d$. 
  If $I = (a_1, \ldots, a_n)$ for appropriate elements
  $a_1, \ldots, a_n \in A$ then it is 
  \begin{displaymath}
    E[I] = \Kern ( E \rTo^{e_{a_1}, \ldots, e_{a_n}} E \times_S \cdots
    \times_S E).
  \end{displaymath}
  In the affine case $S = \spec R$ we have
  \begin{displaymath}
    E[I] = \spec R[X]/(e_{a_1}(X), \ldots,e_{a_n}(X)). 
  \end{displaymath}
\item If $I,J$ are coprime ideals in $A$, then
  \begin{displaymath}
    E[IJ] \cong E[I] \times_S E[J].
  \end{displaymath}
\item If $I$ is coprime to the characteristic of the Drinfeld module
  $E$ then $E[I]$ is \'etale over $S$.
\item The group scheme $E[I]$ is compatible with base change, that is for
  each scheme $T/S$ we have
  \begin{displaymath}
    E[I] \times_S T \cong (E \times_S T)[I].
  \end{displaymath}
  \end{enumerate}
\end{Properties}

\begin{Beweis}
  Cf. \cite{lehmkuhl}, chapter 2, proposition 4.1, page 27 et seq.
\end{Beweis}

If $S=\spec R$ is an affine Scheme and if $\Garbe{L}$ ist trivial we can use the following lemma to describe the group scheme $E[I]$ by a unique additive polynomial $h_I(X)$ of degree $|A/I|^d$.

\begin{Lemma}\label{5.6}
  Let $R$ be an $\Fq$-algebra.
  Let $H \subset \Gas{R}$ be a finite flat subgroup scheme
  of rank $n$ over $R$. Then there is a uniquely defined normalized additive  
  polynomial  $h \in R[X]$ of degree $n$ such that $H = V(h)$.
\end{Lemma}

\begin{Beweis}
  Cf. \cite{lehmkuhl}, chapter 1, lemma 3.3, page 9. 
\end{Beweis}

If $S = \spec L$ is a field then
the characteristic of the Drinfeld module $(L,e)$ is a prime ideal of $A$. Thus we can define the \emph{height} of
$(L,e)$, denoted by $h$.

In the case of an algebraically closed field we have the following explicit description of
the $I$-division points:

\begin{Satz}\label{5.3}
  Let $0 \neq \prim{p} \subset A$ be a prime ideal and let $I =
  \prim{p}^n$ for an $n > 0$. Then we have
  \begin{displaymath}
    E[\prim{p}^n](L) \cong
    \begin{cases}
      (\prim{p}^{-n}/A)^d  & \mathrm{for}\; \prim{p} \neq  \Char E\\
      (\prim{p}^{-n}/A)^{d - h} &\mathrm{for}\; \prim{p} = \Char E.
    \end{cases}
  \end{displaymath}
\end{Satz}

\begin{Beweis}
  Cf. \cite{lehmkuhl}, chapter 2, corollary 2.4, page 24.
\end{Beweis}

This result motivates the following definition:
\begin{Definition}[\cite{drinfeld76}]\label{Def:A}
  Let $E = (\Gas{\Garbe{L}}, e)$ be a Drinfeld module of rank $d$ 
  over $S$ and let $0 \neq I \subsetneq A$ be an ideal. 
  A Level $I$ structure is an $A$-linear map
  \begin{eqnarray*}
    \iota: (I^{-1}/A)^d \rTo E(S),
  \end{eqnarray*}
  such that for all prime ideals $\prim{p} \supseteq I$ we have an identity of Cartier divisors
  \begin{eqnarray*}
    E[\prim{p}] = \sum_{x \in (\prim{p}^{-1}/A)^d} \iota(x).
  \end{eqnarray*}
\end{Definition}

\begin{Remark}
 \begin{enumerate}
  \item If $I$ is coprime to the characteristic of $E$ then a level $I$ structure is an isomorphism of group schemes
  \[
    (I^{-1}/A)^d_S \simeq E[I].
  \]
  \item If $E = (R,e)$  the equality of Cartier divisors simply means an equaltiy of the polynomials
  \[
   h_{\prim{p}}(X) = \prod_{x \in (\prim{p}^{-1}/A)^d} (X-\iota(x)).
  \]
\end{enumerate}
\end{Remark}

\section{Level Structures and Deformations}
In \cite{lehmkuhl}, 3, prop. 3.3 we have the following result.

\begin{Proposition}
 If $(E,\iota)$ is a Drinfeld module equipped with a level $I$ structure $\iota$ then we have the identity of
  Cartier divisors
  \[
    E[I] = \sum_{x \in (I^{-1}/A)^d} \iota(x).
  \]
\end{Proposition}

The proof in loc.~cit. is based on the construction of deformation spaces of Drinfeld modules, isogenies and 
level $I$ structures. In the following we will repeat the basic definitions and results.

\begin{Definition}
\begin{enumerate}
\item Let $i : A \rTo O$ be a complete noetherian $A$-algebra with residue field $\ell$. Let $\mathcal{C}_O$ be the
  category of 
 local artinian $O$-algebras with residue field $\ell$ and let $\hat{\mathcal{C}}_O$ be the category of noetherian complete
 local $O$-algebras with residue field $\ell$.
 \item Let $E_0$ be a Drinfeld module of rank $d$ over $\ell$, and let $B$ be an algebra in $\mathcal{C}_O$. A
  \emph{deformation} of $E_0$  over $B$ is a Drinfeld module of rank $d$ over $\spec B$ which specializes mod
   $\mathfrak{m}_B$ to $E_0$. Thus we obtain a functor:
   \begin{equationarray*}{rcl}
    \mathrm{Def}_{E_0} : \mathcal{C}_O &\rTo& \mathrm{Sets}\\
    B &\rMapsto& \{ \mathrm{Isomorphyclasses}\; \mathrm{of}\; \mathrm{Deformations}\; \mathrm{of}\; E_0\}
   \end{equationarray*}
 \item Let $\varphi_0 : E_0 \rTo F_0$ be an isogeny of Drinfeld modules of rank $d$ over $\ell$. A \emph{deformation} of 
 $\varphi_0$ over $B$ is an isogenie $\varphi: E \rTo F$ where $E$, $F$ are deformations of $E_0$ and $F_0$, such that
 $\varphi$ specializes mod $\mathfrak{m}_B$ to $\varphi_0$. We obtain a corresponding functor:
  \begin{equationarray*}{rcl}
    \mathrm{Def}_{\varphi_0} : \mathcal{C}_O &\rTo& \mathrm{Sets}\\
    B &\rMapsto& \{ \mathrm{Isomorphyclasses}\; \mathrm{of}\; \mathrm{Deformations}\; \mathrm{of}\; \varphi_0\}
    \end{equationarray*}
 \item Let $(E_0, \iota_0)$ be a Drinfeld module of rank $d$ over $\ell$ equipped with a level $I$ structure $\iota_0$. A
 \emph{deformation} is a Drinfeld module $(E, \iota)$ over $B$ of rank $d$ equipped with an level $I$ structure $\iota$ 
 such that $E$ is a deformation of $E_0$ and $\iota$ specializes to $\iota_0$ mod $\mathfrak{m}_B$. We define the functor:
  \begin{equationarray*}{rcl}
    \mathrm{Def}_{(E_0,\iota_0)} : \mathcal{C}_O &\rTo& \mathrm{Sets}\\
    B &\rMapsto& \{ \mathrm{Isomorphyclasses}\; \mathrm{of}\; \mathrm{Deformations}\; \mathrm{of}\; (E_0,\iota_0)\}
    \end{equationarray*}
    \item  We will denote the tangent spaces of the functors above by 
    \[
     T_{E_0} :=  \mathrm{Def}_{E_0}(\ell[\varepsilon]),\; 
     T_{\varphi_0} :=  \mathrm{Def}_{\varphi_0}(\ell[\varepsilon]),\; 
     T_{(E_0,\iota_0)} :=  \mathrm{Def}_{(E_0,\iota_0)}(\ell[\varepsilon])
    \]
where $\ell[\varepsilon]$ is the $\ell$-algebra with $\varepsilon^2 = 0$.
\end{enumerate}
\end{Definition}

\begin{Results}
 \begin{enumerate}
  \item The deformation functor $\mathrm{Def}_{E_0}$ is pro-represented by the smooth $O$-algebra 
  $R_0 := O[[T_1,\ldots T_{d-1}]]$.
  \item The deformation functor $\mathrm{Def}_{\varphi_0}$ is pro-represented by an object in $\hat{\mathcal{C}}_O$.
  \item The deformation functor $\mathrm{Def}_{(E_0,\iota_0)}$ is pro-represented by an object in $\hat{\mathcal{C}}_O$.
 \end{enumerate}

\end{Results}

\section{Equivalence of Definitions}
The question is, what happens if we would change the definition \ref{Def:A} of level $I$ structures with:

\begin{Definition}\label{Def:B}
  Let $E = (\Gas{\Garbe{L}}, e)$ be a Drinfeld module of rank $d$ 
  over $S$ and let $0 \neq I \subsetneq A$ be an ideal. 
  A Level $I$ structure is an $A$ linear map
  \begin{eqnarray*}
    \iota: (I^{-1}/A)^d \rTo E(S),
  \end{eqnarray*}
  such that we have an identity of Cartier divisors
  \begin{eqnarray*}
    E[I] = \sum_{x \in (I^{-1}/A)^d} \iota(x).
  \end{eqnarray*}
\end{Definition}

To distinguish the two definitions we will refer the original definition as \textbf{A} and the one above as \textbf{B}.

\begin{Proposition}\label{prop:1}
Definition \textbf{A} and definition \textbf{B} are equivalent.
\end{Proposition}

As it is proved in \cite{lehmkuhl}, 3, prop. 3.3. being a level $I$ structure in the sense of definition \textbf{A} implies being on in the sense of definition \textbf{B}.
One the other hand it is proved in \cite[Ch.~6, Prop.~6.7]{wiedmanndiss} by a simple counting argument that \textbf{B} implies \textbf{A} if the base scheme $S$ is reduced. On the other hand the result is clear if $I$ is away from the characteristic of the Drinfeld module. Thus for the proof of proposition \ref{prop:1} we are allowed to make the following assumptions:

\begin{enumerate}
 \item $S = \spec B$ where $B$ is the localization of a finitely generated $A$-algebra at a maximal prime Ideal, 
 $\mathfrak{p} := A \cap \mathfrak{m}_B$ is not zero and $\ell := B/\mathfrak{m}_B$ is a finite extension of 
 $A/\mathfrak{p}$.
 \item The result is true if it is true for all quotients $B/\mathfrak{m}_B^n$. So we can assume, that $B$ is a local 
 artinian ring with residue field $\ell$.
 \item We will fix an element $\varpi_\prim{p} \in A$, such that $(\varpi_\prim{p}) = \prim{p}J$ with 
 $\prim{p} \nmid J$. Then we have 
 \[
  A \rInto \hat{A}_\prim{p} \cong A/\prim{p}[[\varpi_\prim{p}]] \rInto \ell[[\varpi_\prim{p}]] =:O
 \]
 such that $\mathfrak{m}_O \cap A = \prim{p}$ and $\mathfrak{p}O = \mathfrak{m}_O$.
 \item As $B$ is artinian, and therefore complete, there is a unique lift of the coefficient field $\ell$ to $B$ and we can 
 consider $B$ as an object of $\mathcal{C}_O$.
 \item We fix a Drinfeld module $E_0 = (e^{(0)},\ell)$ and a level $\prim{p}^n$ structure $\iota_0$. As $\ell$ is a field
 there is no difference between the definitions \textbf{A} and \textbf{B}.
 \item Let $E$ be the universal deformation of $E_0$ and $O$ as above. Then $E$ is defined over
 \[
  R_0 := O[[T_1,\ldots,T_{d-1}]] \cong \ell[[T_0, \ldots, T_{d-1}]]
 \]
for $T_0:=\varpi_\prim{p}$. It is a complete regular local ring of dimension $d$. In especially $R_0$ is integral and the map 
$\Char: A \rTo R_0$ is injective and the Drinfeld module has general characteristic.
\end{enumerate}

To prove the proposition we will follow the arguments of \cite[3.3.1]{lehmkuhl}. The main difference is now to use $\prim{p}^n$ instead of $\prim{p}$.

Let $ p_\prim{p^n}: E \rTo E/E[\prim{p}^n]$ the canonical quotient isogeny of Drinfeld modules with kernel $ E[\prim{p}^n]$. 
The corresponding polynomial $p_\prim{p^n}^\sharp = h_\prim{p^n} \in R_0[X]$ is an additive, normalized and separable polynomial of degree $|A/\prim{p}^n|^d$. We define $L$ to be a splitting field of $h_\prim{p^n}$ over the field of quotients
$\Quot(R_0)$. Using the zeros $V(h_\prim{p^n})$ of $h_\prim{p^n}$ in $L$ we define the $R_0$-algebra 
$R_{h_\prim{p^n}} := R_0[V(h_\prim{p^n})]$ inside $L$. It is an integral and finite extension of $R_0$ because $h_\prim{p^n}$ is normalized over $R_0$.

Along the lines of \cite{lehmkuhl} we will prove

$$ R_{h_\prim{p^n}} \cong R_n $$ 

where $R_n$ is the base ring of the universal Drinfeld module of deformations of $E_0$ and a level $\prim{p}^n$ structure over $\ell$. We will use induction on $n$:
 
If $n=0,1$ then nothing is to prove. For $n>1$ we have:
\[
 R_{n} := R_{n-1}[[S_1,\ldots,S_d]]/\mathfrak{a}
\]
where $\mathfrak{a}$ is the ideal generated by elements of the form 
$e^\sharp_{\varpi_{\prim{p}}}(S_i) - \iota_r(x_i) +  e^\sharp_{\varpi_{\prim{p}}}(\tilde{y}_i)$. These elements are normalized polynomials in $S_i$, so 
\[
 R_{n} := R_{n-1}[S_1,\ldots,S_d]/\mathfrak{a}
\]

We can find elements $\tilde{x}_1,\ldots, \tilde{x}_d$ in $R_{h_\prim{p^n}}$ such that $e^\sharp_{\varpi_{\prim{p}}}(\tilde{x}_i) = x_i$ and $\tilde{x}_1,\ldots, \tilde{x}_d$ generates $V(h_\prim{p^n})$ as an $A/\prim{p}^n$-module. We define a map of $R_{n-1}$-algebras:
\begin{equationarray*}{rcl}
 R_{n}&\rTo&R_{h_\prim{p^n}}\\
 S_i &\rMapsto& \tilde{x}_i - \tilde{y_i}
\end{equationarray*}
By assumption the map is well defined and surjective. As both rings have dimension $d$, it is an isomorphism.

\begin{Corollary}
 $R_{h_\prim{p^n}}$ is a regular algebra in $\hat{\mathcal{C}}_O$.
\end{Corollary}

Now we can use the setup of \cite{lehmkuhl}, proof of proposition 3.3.1, to show that 
$E_{\prim{p}^n} := E \otimes_{R_0} R_{h_\prim{p^n}}$ and a corresponding lift of $\iota_0$ is the universal deformation of level $\prim{p}^n$. This is true if we can prove \cite{lehmkuhl}, lemma 3.3.2, with $\prim{p}$ replaced by $\prim{p}^n$, but there is no obstruction to do so.

Now we are done, because $R_{h_\prim{p^n}}$ is an integral ring and definition $A$ and $B$ coincide on the level of the universal deformation.



\begin{thebibliography}{Wie04}

\bibitem[Dri76]{drinfeld76}
V.~G. Drinfeld.
\newblock Elliptic modules.
\newblock {\em Math. USSR Sbornic}, 23:561--592, 1976.

\bibitem[KM85]{KatzMazur}
Nicholas~M. Katz and Barry Mazur.
\newblock {\em Arithmetic moduli of elliptic curves}, volume 108 of {\em Annals
  of Mathematics Studies}.
\newblock Princeton University Press, Princeton, NJ, 1985.

\bibitem[Leh09]{lehmkuhl}
Thomas Lehmkuhl.
\newblock Compactification of the {D}rinfeld modular surfaces.
\newblock {\em Mem. Amer. Math. Soc.}, 197(921), 2009.

\bibitem[Wie04]{wiedmanndiss}
Stefan Wiedmann.
\newblock {\em {Drinfeld modules and elliptic sheaves. (Drinfeld-Moduln und
  elliptische Garben.)}}.
\newblock PhD thesis, {G\"ottingen: Univ. G\"ottingen,
  Mathematisch-Naturwissenschaftliche Fakult\"aten (Dissertation). 98~p. },
  2004.

\end{thebibliography}
\end{document}